\documentclass[a4, 11pt] {amsart} 
\usepackage{amssymb,times, amscd,amsmath,amsthm, xypic}
\usepackage{geometry}

\author{Jean-Paul Brasselet}

\address{Institute de Math\'{e}matiques de Marseille,
UMR 7373\\
Campus de Luminy - Case 907, 13288 Marseille Cedex 9, France}
\email{jean-paul.brasselet@univ-amu.fr}

\author{J\"{o}rg Sch\"{u}rmann$^{(*)}$}

\address{Westf\"{a}lische Wilhelms-Universit\"{a}t,
Mathematische Institut\\
Einsteinstrasse 62,
48149 M\"{u}nster, Germany}
\email{jschuerm@math.uni-muenster.de}

\author{Shoji Yokura$^{(**)}$}

\address{Department of Mathematics and Computer Science, 
Graduate School of Science and Engineering, Kagoshima University,
1-21-35 Korimoto, Kagoshima 890-0065, Japan}
\email{yokura@sci.kagoshima-u.ac.jp}

\date{}
\thanks {Mathematics Subject Classification: 14C17, 14C40, 14F25, 14F45, 14Q15, 32S35}
\thanks{(*) Partially supported by the SFB 878 groups, geometry and actions}
\thanks{(**) Partially supported by JSPS KAKENHI Grant Numbers 23244008, 24340007, 24540085 }

\title 
{Motivic and 
derived motivic Hirzebruch classes}
\begin{document} 
\numberwithin{equation}{section}
\newtheorem{thm}[equation]{Theorem}
\newtheorem{pro}[equation]{Proposition}
\newtheorem{prob}[equation]{Problem}
\newtheorem{qu}[equation]{Question}
\newtheorem{cor}[equation]{Corollary}
\newtheorem{con}[equation]{Conjecture}
\newtheorem{lem}[equation]{Lemma}
\theoremstyle{definition}
\newtheorem{ex}[equation]{Example}
\newtheorem{defn}[equation]{Definition}
\newtheorem{ob}[equation]{Observation}
\newtheorem{rem}[equation]{Remark}
\renewcommand{\rmdefault}{ptm}
\def\alp{\alpha}
\def\be{\beta}
\def\jeden{1\hskip-3.5pt1}
\def\om{\omega}
\def\bigstar{\mathbf{\star}}
\def\ep{\epsilon}
\def\vep{\varepsilon}
\def\Om{\Omega}
\def\la{\lambda}
\def\La{\Lambda}
\def\si{\sigma}
\def\Si{\Sigma}
\def\Cal{\mathcal}
\def\m {\mathcal}
\def\ga{\gamma}
\def\Ga{\Gamma}
\def\de{\delta}
\def\De{\Delta}
\def\bF{\mathbb{F}}
\def\bH{\mathbb H}
\def\bPH{\mathbb {PH}}
\def \bB{\mathbb B}
\def \bA{\mathbb A}
\def \bC{\mathbb C}
\def \bOB{\mathbb {OB}}
\def \bM{\mathbb M}
\def \bOM{\mathbb {OM}}
\def \mA{\mathcal A}
\def \mB{\mathcal B}
\def \mC{\mathcal C}
\def \mR{\mathcal R}
\def \mH{\mathcal H}
\def \mM{\mathcal M}
\def \mM{\mathcal M}
\def \mT{\mathcal {T}}
\def \mAB{\mathcal {AB}}
\def \bK{\mathbb K}
\def \bG{\mathbf G}
\def \bL{\mathbb L}
\def\bN{\mathbb N}
\def\bR{\mathbb R}
\def\bP{\mathbb P}
\def\bZ{\mathbb Z}
\def\bC{\mathbb  C}
\def \bQ{\mathbb Q}
\newcommand{\bb}[1]{\mbox{$\mathbb{#1}$}}

\def\op{\operatorname}

\maketitle

\begin{abstract}
In this paper we give a formula for the Hirzebruch $\chi_y$-genus $\chi_y(X)$ and similarly for the motivic Hirzebruch class $T_{y*}(X)$ for possibly singular varieties $X$, using the Vandermonde matrix. 
Motivated by the notion of secondary Euler characteristic and higher Euler characteristic, we consider a similar notion for the motivic Hirzebruch class, which we call a \emph{derived motivic Hirzebruch class}.
\end{abstract}

\tableofcontents

\section{Introduction}
First we will recall that the Euler-Poincar\'e characteristic is a kind of ``generalization" or ``broad extension" of the counting of a finite set, where 
the counting of a finite set $X$ is the so-called cardinality, i.e.
$$c(X):= |X| = \, \, \text{the number of the elements in the set $X$}.$$

Certainly the counting  $c$ for finite sets satisfies the following basic properties:

\begin{enumerate}

\item
$\displaystyle  A \cong A'$ (bijection or equipotent) $\Longrightarrow$ $c(A) = c(A')$,

\item
$c(A) = c(A \setminus B) + c(B)$ for  $B \subset A$,  (this is called \emph{``scissor formula"} or \emph{``motivic"})

\item
$c(A \times B) = c(A)\cdot c(B)$,

\item
$c(pt) = 1$. \quad (Here $pt$ denotes one point.) 
\end{enumerate}

Now, if we consider the following ``topological counting" $c$ on the category of some ``nice" topological spaces such that $c(X) \in \bZ $ and it satisfies the following four properties:
\begin{itemize}
\item
$ X\cong X'$ (homeomorphism = $\Cal{TOP}$- isomorphism) $\Longrightarrow$ $c(X) = c(X')$,

\item
$c(X) = c(X\setminus Y) + c(Y)$ for $Y \subset X$

\item
$c(X \times Y) = c(X)\cdot c(Y)$,

\item
$c(pt) = 1,$ 
\end{itemize}
then one can show that if such a $c$ exists, then we must have that
$$c(\bR^1) =  -1, \quad \text {hence} \quad c({\bR}^n) = (-1)^n.$$
Hence,  if $X$ is a finite $CW$-complex with $\sigma_n(X)$ denoting the number of open $n$-cells, then
$$c(X) = \sum_n (-1)^n \sigma_n(X) = \chi(X)$$
is the Euler--Poincar\'e characteristic of $X$. Namely, the topological counting $c$ is uniquely determined and it is the compactly supported Euler--Poincar\'e characteristic.

\begin{rem} 
\begin{enumerate}
\item Such a counting is not defined for all topological spaces, as one can see for example that such a $c$ is not well-defined on the discrete space $\mathbb Z$ of integers. Such a counting is defined on ``nice" spaces such as finite $CW$-complexes. Here we do not bother ourselves to specify what we mean by ``nice" (cf. C. Peters' TATA Lecture Notes \cite{Peters}).
\item It would be safe to say that the reason why the Euler-Poincar\'e characteristic (which is the very fundamental, basic but still important topological invariant in topology, geometry and physics) is defined by \emph{the alternating sum} of the numbers of vertices, edges, faces, and so on,
$$\chi(X) := V - E + F - \cdots.$$
\emph{comes from our usual simple ``counting"}.
\end{enumerate}
\end{rem}

Now, let us consider such a counting on the the category $\m V$ of algebraic varieties:
\begin{itemize}
\item
$ X\cong X'$ ($\m V$-isomorphism) $\Longrightarrow$ $c(X) = c(X')$,

\item
$c(X) = c(X\setminus Y) + c(Y)$ for a closed subvariety $Y \subset X$

\item
$c(X \times Y) = c(X)\cdot c(Y)$,

\item
$c(pt) = 1.$ \\
\end{itemize}
If such an ``algebraic" counting $c$ exists, then it follows from the 
decomposition of the $n$-dimensional complex projective space
$$\bP^n = \bC^0\sqcup \bC^1 \sqcup \cdots \sqcup \bC^{n-1} \sqcup \bC^n$$
that we must have 
$$c(\bP^n) = 1 - y + y^2 - y^3 + \cdots + (-y)^n$$
where $y := -c(\bC^1) \in \bZ$.
In fact, it follows from Deligne's mixed Hodge structure that the following Hodge--Deligne polynomial

$$\chi _{u,v}(X) := \sum_{i, p, q \geq 0} (-1)^i (-1)^{p+q}\op{dim}_{\bC} (Gr^p_F Gr^W_{p+q} H^i_c(X, \bC)) u^p v^q$$
satisfies the above four properties, namely any Hodge--Deligne polynomial $\chi _{u,v}$ with $uv = -y$ is such a $c$.  The Hirzebruch $\chi_y$ characteristic is nothing but $\chi_{y, -1}$ and the most important and interesting ones are the following:
\begin{itemize}
\item $ y = -1: \quad \text {$\chi_{-1} = \chi$, the topological Euler--Poincar\'e characteristic}$
\item $ y = 0: \quad \text {$\chi_0 = \chi^a$, the arithmetic genus}$ (for a compact nonsingular variety)
\item $ y = 1: \quad \text {$\chi_{1} = \sigma$, the signature}$ (for a compact nonsingular variety) \\
\end{itemize}

\begin{rem} Here we note that $\chi_0=\chi^a$ also holds for $X$ compact with at most Du Bois singularities (by \cite{BSY1}) and $\chi_1=\sigma$ for $X$ a projective rational homology manifold (by \cite[\S 3.6]{MSS}). 
\end{rem}
It turns out (see \cite{BSY1} and \cite{YokuraMSRI}) that  the Hodge--Deligne polynomial $\chi_{u, v}: K_0(\m V) \to \bZ[u,v]$ can be extended as a class version only when $u=y, v =-1$, just like Hirzebruch--Riemann--Roch was extended by A. Grothendieck as a natural transformation from the covariant functor of coherent sheaves to the rational homology theory, which is called Grothendieck--Riemann--Roch. Here $K_0(\m V)$ is the Grothendieck group of complex algebraic varieties with respect to the scissor relation. 
Namely only the Hirzebruch $\chi_y$ characteristic
$$\chi_y: K_0(\m V) \to \bZ[y]$$
can be extended as a class version
$$T_{y*}: K_0(\m V/X) \to H_*(X)\otimes \bQ[y].$$
This is called the \emph{motivic Hirzebruch class}. Here $K_0(\m V/X)$ is the relative Grothendieck group of complex algebraic varieties as recalled in \S 3 and $H_*(X)$ is the Borel-Moore homology group.

In this paper we give some formulas for the motivic Hirzebruch class and its ``derived version", which we call ``derived motivic Hirzebruch class", motivated by \emph{higher Euler characteristic} generalizing the \emph{secondary characteristic} (see \cite{Ramachandran} and cf.\cite{GN}).

\section{Hirzebruch $\chi_y$-genus and Hirzebruch class $T_y$}

First we recall the  definition of the Hirzebruch $\chi_y$-genus. Let $X$ be a smooth complex projective
variety. The
\emph{$\chi_{y}$-genus} of $X$ is defined by 
$$\chi_{y}(X):= 
\sum_{p\geq 0} \chi(X,\Lambda^{p}T^{*}X)y^p
= \sum_{p\geq 0} \left( \sum_{i\geq 0}
(-1)^{i}\op{dim}_{\mathbb C}H^{i}(X,\Lambda^{p}T^{*}X) \right)y^p\:,
$$

Thus the $\chi_y$-genus is the generating function of the Euler-Poincar\'e characteristic $\chi(X,\Lambda^{p}T^{*}X)$ of 
the $p$-th exterior power $\Lambda^{p}T^{*}X$ of the cotangent bundle $T^*X$, which shall be simply denoted by $\chi^p(X)$:
$$\chi_{y}(X)= \sum_{p\geq 0} \chi^p(X)y^p.$$
Since $\Lambda^{p}T^{*}X =0$ for $p> \op{dim}_{\mathbb C}X$, $\chi_{y}(X)$ is a polynomial of at most degree $\op{dim}_{\mathbb C} X$.


More generally, for $E$ a holomorphic vector bundle over $X$, the Hirzebruch $\chi_y$-genus of $E$ is defined by
$$\chi_{y}(X,E):= 
\sum_{p\geq 0} \chi(X,E\otimes \Lambda^{p}T^{*}X) y^{p}\\
=\sum_{p\geq 0} \left( \sum_{i\geq 0}
(-1)^{i}\op{dim}_{\mathbb C}H^{i}(X,E\otimes \Lambda^{p}T^{*}X) \right) y^{p} \:.
$$
Then we have 
\begin{equation} \label{eq:gHRR} 
\chi_{y}(X,E)= \int_{X} T_{y}(TX)\cdot ch_{(1+y)}(E) \cap [X]
\quad \in \bb{Q}[y], \tag{gHRR}
\end{equation}
\[\text{with} \quad ch_{(1+y)}(E):= \sum_{j=1}^{rank\;E} e^{\beta_{j}(1+y)}
\quad \text{and} \quad T_{y}(TX):= \prod_{i=1}^{dim X}
Q_{y}(\alpha_{i})  \:.\]
Here $\beta_{j}$ are the Chern roots of $E$, $\alpha_{i}$ are the Chern
roots of the tangent bundle $TX$, and $Q_{y}(\alpha)$ is the normalized
power series 
\[Q_{y}(\alpha):= \frac{\alpha(1+y)}{1-e^{-\alpha(1+y)}} -\alpha y
\quad \in \bb{Q}[y][[\alpha]] \:. \]
Note that this power series $Q_{y}(\alpha)$ specializes to
\begin{displaymath}
Q_{y}(\alpha) = 
\begin{cases}
\:1+\alpha &\text{for $y=-1$,}\\
\:\frac{\alpha}{1-e^{-\alpha}} &\text{for $y=0$,}\\
\:\frac{\alpha}{\tanh \alpha} &\text{for $y=1$.}
\end{cases} \end{displaymath}
Therefore the {\em modified Todd class\/} $T_{y}(TX)$ unifies the
following important characteristic cohomology classes of $TX$:
\begin{displaymath}
T_{y}(TX) = 
\begin{cases}
\:c(TX) &\text{the total {\em Chern class\/} for $y=-1$,}\\
\:td(TX) &\text{the total {\em Todd class\/} for $y=0$,}\\
\:L(TX) &\text{the total {\em Thom-Hirzebruch L-class\/} for $y=1$.}
\end{cases} \end{displaymath}
We call the modified Todd class $T_{y}(TX)$ \emph{the Hirzebruch class of $X$}.

The coefficient of the power $y^p$ of the Hirzebruch class $T_y(E)$ shall be denoted by $T^p(E)$ (cf. Hirzebruch's book \cite{Hirzebruch}):
$$ T_y(E)= \sum_{i=0}^{\op{rank} E} T^p(E)y^p.$$
Here we emphasize that each $T^p(E)$ is a polynomial of Chern classes or Pontryagin classes without the variable $y$ not involved at all. 

The Hirzebruch $\chi_y$-genus is by the generalized Hirzebruch--Riemann--Roch formula given by
$$\chi_y(X) = \int_X T_y(TX) \cap [X].$$
Hence for a compact nonsingular variety $X$ of dimension $n$
$$\chi_y(X) = \sum_{i=0}^n \left (\int_X T^p(TX) \cap [X] \right ) y^p.$$
So we note that
$$\chi^p(X) =  \int_X T^p(TX) \cap [X].$$
Since we eventually deal with homology classes, we define 
$$T^p_*(X):= T^p(TX) \cap [X], \qquad T_{y*}(X) := \sum_{p=0}^n T_*^p(X)y^p.$$

For the distinguished three values $-1, 0, 1$ of $y$, by the definition we have the following:

\begin{itemize}
\item $c(E) = T_{-1}(E) = T^0(E) - T^1(E) + T^2(E) - \cdots +(-1)^n T^{\op{rank} E}(E),$

\item $td(E) = T_0(E) = T^0(E), \hspace{3cm}$

\item $L(E) = T_1(E) = T^0(E) + T^1(E) + T^2(E) + \cdots + T^{\op{rank} E}(E).$
\item $\chi(X) = \chi_{-1}(X) = \chi^0(X)- \chi^1(X) + \chi^2(X)- \cdots +(-1)^n\chi^n(X),$

\item $\chi^a(X) = \chi_0(TX) = \chi^0(X), \hspace{3cm}$

\item $\sigma(X) = \chi_1(TX) = \chi^0(X)+ \chi^1(X) + \chi^2(X)+\cdots +\chi^n(X).$
\end{itemize}

\begin{rem}\label{inversion} 
\begin{enumerate}
\item Each individual coefficient $T^p(E)$ of the Hirzebruch class $T_y(E)$ may be computed from the very definition of it, but as the above formulas for these three distinguished values suggest, certain summations of all these coefficients make more sense than each individual $T^p(E)$ geometrically or topologically. It is the same for the Hirzebruch $\chi_y$-genus $\chi_y(X)$. 

\item It would be worthwhile to observe that 
$$\frac{L(E) + c(E)}{2} = T^0(E) + T^2(E) + T^4(E) + \cdots: \quad \text{the ``even part"}$$
$$\frac{L(E) - c(E)}{2} = T^1(E) + T^3(E) + T^5(E) + \cdots: \quad \text{the ``odd part"}.$$
\item It follows from \cite{Hirzebruch} that for a compact nonsingular variety $X$
$$\chi^p(X) = (-1)^n \chi^{n-p}(X).$$
It is also known (e.g., see \cite{Kotschick2} and cf. \cite{Kotschick1}) that $\chi^0, \chi^1, \cdots, \chi^{[\frac{n}{2}]}$ are linearly independent, which means that for any compact nonsingular variety $X$ of dimension $n$
$$r_0\chi^0(X) + r_1 \chi^1(X) +  \cdots + r_{[\frac{n}{2}]} \chi^{[\frac{n}{2}]}(X) = 0$$
implies that $r_0 = r_1 = \cdots = r_{[\frac{n}{2}]} =0$.
As a corollary of this, we can also say that $T_*^0, T_*^1, \cdots, T_*^{[\frac{n}{2}]}$ are linearly independent, otherwise that the linear dependence of  $T_*^0, T_*^1, \cdots, T_*^{[\frac{n}{2}]}$ implies the linear dependence of $\chi^0, \chi^1, \cdots, \chi^{[\frac{n}{2}]}$, which contradicts to the above linear independence.

\item 
The above ``duality formula" $\chi^p(X) = (-1)^n \chi^{n-p}(X)$ implies the following ``inversion formula": for a compact nonsingular variety $X$ 
$$\chi_y(X) = (-y)^n \chi_{\frac{1}{y}}(X), \quad \text{i.e.,} \quad \chi_{\frac{1}{y}}(X) = \left (-\frac{1}{y} \right )^n \chi_y(X).$$
\end{enumerate}
\end{rem}

\section{Motivic Hirzebruch classes $T_{y*}$}
The Hirzebruch $\chi_y$-genus was extended to the case of singular varieties, using Deligne's mixed Hodge structures, i.e., 
\begin{align*}
\chi_y(X) &:= \sum_{i, p \geq 0} (-1)^i \op{dim}_{\mathbb C} Gr^p_{\mathcal F} \left( H^i_c (X, \mathbb C) \right) (-y)^p\\
& = \sum_{p \geq 0} \left (\sum_{i \geq 0} (-1)^{i+p} \op{dim}_{\mathbb C} 
Gr^p_{\mathcal F} \left ( H^i_c (X, \mathbb C) \right )\right)  y^p.
\end{align*}
Here $\mathcal  F$ is the Hodge filtration in the mixed Hodge structure of  $H^i_c(X, \mathbb C) $. Thus for a possibly singular variety $X$, the coefficient $\chi^p(X)$ of the above Hirzebruch $\chi_y$-genus  $\chi_y(X)$ is
$$\chi^p(X) = \sum_{i \geq 0} (-1)^{i+p} \op{dim}_{\mathbb C} 
Gr^p_{\mathcal F} \left ( H^i_c(X, \mathbb C) \right ).$$
Here we remark that the degree of the above integral polynomial $\chi_y(X)$ of a possibly singular variety is also at most the dimension of $X$ like in the smooth case (cf. \cite[Corollary 3.1 (1)]{BSY1}).

The three distinguished characteristic classes have been extended to the case of possibly singular varieties as natural transformations from certain covariant functors to the homology functor. This formulation is analogous to the interpretation that the classical theory of characteristic classes of vector bundles is a natural transformation from the contravariant monoid functor $\m Vect$ or the Grothendieck $K$-theory of real or complex vector bundles to the contravariant cohomology theory. 
The three theories of characteristic classes of singular varieties are the following:
\begin{enumerate}
\item MacPherson's Chern class transformation \cite{MacPherson}:
$$c_*: F(X) \to H_*(X),$$ 
where $F$ is the covariant functor assigning to $X$ the abelian group $F(X)$ of constructible functions on $X$. Here we remark that J.-P. Brasselet and M.-H. Schwartz \cite{BS} (see also \cite{AB})
showed that MacPherson's Chern class $c_*(\jeden_X)$ corresponds to the Schwartz class $c^S(X) \in H^{*}_X(M) = H^*(M, M \setminus X)$ (see \cite{Schw1, Schw2}) by  Alexander duality for $X$ embedded in the smooth complex manifold $M$. That is why the total homology class $c_*(X) := c_*(\jeden_X)$ is called the Chern--Schwartz--MacPherson class of $X$.

\item Baum--Fulton--MacPherson's Todd class or Riemann--Roch \cite{BFM}:
$$td_*: G_0(X) \to H_*(X)\otimes \bQ,$$
where $G_0$ is the covariant functor assigning to $X$ the Grothendieck group $G_0(X)$ of coherent sheaves on $X$.
\item Goresky-- MacPherson's homology $L$-class \cite{GM}, which is extended as a natural transformation by Sylvain Cappell and Julius Shaneson \cite{CS} (also see \cite{Yokura-TAMS}):
$$L_*: \Omega(X) \to  H_*(X)\otimes \bQ,$$
where $\Omega$ is the covariant functor assigning to $X$ the cobordism group $\Omega(X)$ of self-dual constructible sheaf complexes on $X$.
\end{enumerate}

In our previous paper \cite{BSY1} (see also \cite {BSY2}, \cite{SY}, \cite{S} and \cite{YokuraMSRI}) we introduced the motivic Hirzebruch class 
$${T_y}_*: K_0(\m V/X) \to H_*(X)\otimes \bQ[y],$$
where $K_0(\m V/X)$ is the relative Grothendieck group of the category $\m V$ of complex algebraic varieties, i.e. the free abelian group generated by the isomorphism classes $[V \xrightarrow {h} X]$ of morphism $h \in hom_{\m V}(V, X)$ modulo the relations
\begin{itemize}
\item $[V_1 \xrightarrow {h_1} X] + [V_2 \xrightarrow {h_2} X] =  [V_1 \sqcup V_2 \xrightarrow {h_1 + h_2} X],$ with $\sqcup$ the disjoint union, and
\item  $[V \xrightarrow {h} X] = [V \setminus W \xrightarrow {h_{V \setminus W}} X] +  [ W \xrightarrow {h_{W}}  X]$ for $W \subset V$  a closed subvariety of $V$. 
\end{itemize}

${T_y}_*: K_0(\m V/X) \to H_*(X)\otimes \bQ[y]$ is the unique natural transformation satisfying the normalization condition that if $X$ is non-singular, then
$${T_y}_*{[X \xrightarrow {\op{id}_X}  X]} = T_y(TX) \cap [X].$$
Here $T_y(TX)$ is the above Hirzebruch class.

Our ${T_y}_*: K_0(\m V/X) \to H_*(X)\otimes \bQ[y]$ ``unifies" the above three characteristic classes $c_*, td_*, L_*$ in the sense that we have the following commutative diagrams:
$$\xymatrix{
& K_0(\Cal V/X)  \ar [dl]_{\epsilon} \ar [dr]^{{T_{-1}}_*} \\
{F(X) } \ar [rr] _{c_*\otimes \mathbb Q}& &  H_*(X)\otimes \bQ.}
$$
$$\xymatrix{
&  K_0(\Cal V/X)  \ar [dl]_{\Gamma} \ar [dr]^{{T_{0}}_*} \\
{G_0(X) } \ar [rr] _{td_*}& &  H_*(X)\otimes \bQ.}
$$
$$\xymatrix{
& K_0(\Cal V/X)  \ar [dl]_{\omega} \ar [dr]^{{T_{1}}_*} \\
{\Omega(X) } \ar [rr] _{L_*}& &  H_*(X)\otimes \bQ.}
$$
This ``unification" could be considered as a positive answer to the following remark which is stated at the very end of MacPherson's survey article \cite{MacPherson2} (which is a paper version of his survey talk about characteristic classes of singular varieties at Brazilian Math. Colloquium in 1973):
\emph{``It remains to be seen whether there is a unified theory of characteristic classes of singular varieties like the classical one outlined above."}\footnote{At that time Goresky--MacPherson's homology $L$-class was not available yet and it was defined only after the theory of Intersection Homology was invented by Mark Goresky and Robert MacPherson in 1980. }

\section{Naive explicit formulae for $\chi_y(X)$ and $T_{y*}(X)$}

\subsection{A natural question} Whenever we have given talks about the above motivic Hirzebruch class and $\chi_y$-genus, commenting or emphasizing that our ${T_y}_*: K_0(\m V/X) \to H_*(X)\otimes \bQ[y]$ unifies the above three characteristic classes $c_*, td_*, L_*$, 
the Hirzebruch class $T_y(E)$ specialize to the well-known three distinguished ones: Chern class, Todd class and $L$-class, thus the Hirzebruch $\chi_y$-genus specializes to Euler-Poincar\'e characteristic, the arithmetic genus and the signature for $y=-1, 0, 1$, respectively,
we always have been asked the following question:
\begin{qu} How about  other values $\chi_y$ and $T_{y*}$ for $y \not =-1, 0, 1$,  say, at other integers?
\end{qu}

A motivation of the present work is trying to answer this very reasonable and natural question. Although we have been unable to give a complete answer, the one we give in this paper would be a reasonable one at the moment , considering the fact that (as far as we know) there is no literature available of concrete or explicit formula for $\chi_y(X)$ for an even smooth variety $X$ and for general integers $y$. The idea of our formula is quite simple, because $\chi_y(X)$ and $T_{y*}(X)$ are both polynomials of $y$ of finite degree. Then such a polynomial can be completely described using the special values at $\op{dim} X+1$ points, using the Interpolation Formula or the Vandermonde matrix. 

\subsection{Interpolation polynomial and Vandermonde}

Here we recall some basic things for the sake of completeness.

Let $f: \mathbb R \to \mathbb R$ be a given function and let $\{a_i\}$ ($0 \leqq i \leqq n$) be mutually distinct points, i.e. $a_i \not = a_j (i \not = j)$. An interpolation polynomial for the function $f$ is determined by the following Lagrange interpolation polynomial:
$$ p(x) = \sum_{i=0}^{n} \left (\prod_{0 \leqq j \leqq n, \\
j \not=  i} \frac{x-a_j}{a_i - a_j}\right )f(a_i).$$

Expressing the above interpolation polynomial in the form of  
$$p(x) = p_0 + p_1 x + p_2x^2 + \cdots + p_n x^n$$
can be done directly by using the Vandermonde matrix.
Indeed we have the following linear equations:
\begin{eqnarray*}
\left\{
\begin{array}{lllllll}
p_0 + p_1a_0 + p_2 a_0^2 + p_3 a_0^3 + \cdots + p_n a_0^n = f(a_0)\\
p_0 + p_1a_1+ p_2 a_1^2 + p_3 a_1^3 + \cdots + p_n a_1^n = f(a_1)\\
p_0 + p_1a_2 + p_2 a_2^2 + p_3 a_2^3 + \cdots + p_n a_2^n = f(a_2)\\
 \hspace{1cm} \cdots \cdots \cdots \cdots \cdots \cdots  \cdots \cdots \cdots \\
  \hspace{1cm} \cdots \cdots \cdots \cdots \cdots \cdots  \cdots \cdots \cdots \\
   \hspace{1cm} \cdots \cdots \cdots \cdots \cdots \cdots  \cdots \cdots \cdots \\
p_0 + p_1a_n + p_2 a_n^2 + p_3 a_n^3 + \cdots + p_n a_n^n = f(a_n),\\
\end{array}
\right.
\end{eqnarray*}
i.e.,

\begin{eqnarray*}
\left( 
\begin{array}{ccccccc}
1  & a_0 & a_0^2 & a_0^3 \cdots \cdots & a_0^n \\
1  & a_1 & a_1^2 & a_1^3 \cdots  \cdots & a_1^n \\
1  & a_2 & a_2^2 & a_2^3 \cdots  \cdots & a_2^n \\

\vdots &&\ddots &&\vdots \\
1  & a_n & a_n^2 & a_n^3 \cdots \cdots & a_n^n \\
\end{array} 
\right)
\left( \begin{array}{c}
p_0 \\
p_1\\
p_2 \\
\vdots\\
p_n\\
\end{array}
\right )
= 
\left( 
\begin{array}{c}
f(a_0)\\
f(a_1)\\
f(a_2)\\
\vdots\\
f(a_n)\\
\end{array}
\right)\\
\end{eqnarray*}

Let $V(a_0,a_1,a_2, \cdots, a_n)$ be the above Vandermonde matrix . The determinant of this Vandermonde is
$$\op{det} \left (V(a_0,a_1,a_2, \cdots, a_n) \right )= \prod _{0 \leqq i<j\leqq n} (a_j - a_i),$$
which is $\not = 0$ since $a_i \not = a_j (i \not = j)$. Therefore the coefficients $p_0, p_1, p_2, \cdots, p_n$ can be determined by the following equation, by computing the inverse of the Vandermonde 
$V(a_0,a_1,a_2, \cdots, a_n)$:

\begin{eqnarray*}
\left( 
\begin{array}{c}
p_0 \\
p_1\\
p_2 \\
\vdots\\
p_n\\
\end{array}
\right )
=
\left( 
\begin{array}{ccccccc}
1  & a_0 & a_0^2 & a_0^3 \cdots \cdots & a_0^n \\
1  & a_1 & a_1^2 & a_1^3 \cdots  \cdots & a_1^n \\
1  & a_2 & a_2^2 & a_2^3 \cdots  \cdots & a_2^n \\

\vdots &&\ddots &&\vdots \\
1  & a_n & a_n^2 & a_n^3 \cdots \cdots & a_n^n \\
\end{array} 
\right)^{\LARGE -1}
\left( 
\begin{array}{c}
f(a_0)\\
f(a_1)\\
f(a_2)\\
\vdots\\
f(a_n)\\
\end{array}
\right)\\
\end{eqnarray*}

\begin{rem} We note that the above Lagrange interpolation formula and the method of using the inverse of the Vandermonde matrix can be applied to any function from $\mathbb R$ ($\mathbb C$, resp.) to any vector space $\mathcal {VEC} $ over $\mathbb R$ ($\mathbb C$, resp.).
\end{rem}

\subsection{Explicit computations}
In this section we give an explicit formula for Hirzebruch's $\chi_y$-genus $\chi_y(X)$ and the Hirzebruch class $T_y(X)$ for a compact nonsingular variety $X$ of dimension $n$. Here we recall that
\begin{itemize}
\item $T_{y*}(X) =T^0_*(X) + T^1_*(X)y + T^2_*(X)y^2 + \cdots + T^n_*(X)y^n \in H_*(X)\otimes \mathbb Q[y],$
\item $\chi_y(X) = \chi^0(X)+ \chi^1(X)y + \chi^2(X)y^2+\cdots +\chi^n(X)y^n \in \mathbb Q[y].$
\end{itemize}
We just deal with $\chi_y(X)$ since it is exactly the same for $T_{y*}(X)$.

\begin{ex}
$n=1$: $\chi_y(X) = \chi^0(X)+ \chi^1(X)y$. Since
$$\chi(X) = \chi^0(X)- \chi^1(X),\quad \chi^a(X) = \chi^0(X),\quad \sigma(X) = \chi^0(X)+ \chi^1(X),$$
we get that  $ \chi^1(X) = \frac{\sigma (X) - \chi(X)}{2}.$ Here we note that the signature $\sigma(X) =0$ by definition, since it is defined to be zero if the real dimension $\op{dim}_{\mathbb R}X \not \equiv 0  \op{mod} 4$.
Hence we have
$$\chi_y(X) = \chi^a(X) - \frac{\chi(X)}{2}y .$$
Here we also note that
$$\chi^a(X) = \frac{\chi(X)}{2}.$$
\end{ex}

\begin{ex}
$n=2$: $\chi_y(X) = \chi^0(X)+ \chi^1(X)y + \chi^2(X)y^2$. Since
$$\chi(X) = \chi^0(X)- \chi^1(X) +\chi^2(X),\quad \chi^a(X) = \chi^0(X),\quad \sigma(X) = \chi^0(X)+ \chi^1(X) +\chi^2(X),$$
we get that 
$ \displaystyle \chi^1(X) = \frac{\sigma (X) - \chi(X)}{2}, \, \chi^2(X) = \frac{\chi(X)+ \sigma(X) - 2\chi^a(X)}{2}.$
 Hence we have
$$\chi_y(X) = \chi^a(X) + \frac{\sigma (X) - \chi(X)}{2}y + \frac{\chi(X)+ \sigma(X) - 2\chi^a(X)}{2}y^2.$$
\end{ex}

Since the Hirzebruch $\chi_y$-genus is multiplicative and the Hirzebruch homology class $T_*$ is also multiplicative, i.e. respectively
$$\chi_y(X \times Y) = \chi_y(X) \cdot \chi_y(Y) \quad \text{and} \quad  T_{y*}(X \times Y) = T_{y*}(X) \times T_{y*}(Y),$$
we obtain the following formulas:
\begin{thm} Let $C_i \, \, (1 \leqq i \leqq s)$ be a compact nonsingular curve and $S_j \, \, (1 \leqq i \leqq t)$ be a compact nonsingular surface. Then we have the following formulae:
\begin{align*}
\chi_y(C_1\times  & C_2 \times \cdots   \times C_s \times S_1 \times S_2 \times \cdots \times S_t) \qquad \\
& = \prod_{i=1}^s \left (\chi^a(C_i) - \frac{\chi(C_i)}{2} y \right ) \cdot  \\
& \qquad \qquad 
\prod_{j=1}^t \left (\chi^a(S_j) + \frac{\sigma (S_j) - \chi(S_j)}{2} y  +  \frac{\sigma (S_j) + \chi(S_j) - 2\chi^a(S_j)}{2} y^2\right ).
\end{align*}
\begin{align*}
T_{y*}(C_1\times & C_2 \times \cdots   \times C_s \times S_1 \times S_2 \times \cdots \times S_t) \\
& = \prod_{i=1}^s \left (Td_*(C_i) + \frac{L_*(C_i) - c_*(C_i)}{2} y \right ) \times  \\
& \qquad \qquad  \prod_{j=1}^t \left (Td_*(S_j) + \frac{L_* (S_j) - c_*(S_j)}{2} y  +  \frac{L_* (S_j) + c_*(S_j) - 2Td_*(S_j)}{2} y^2\right ).
\end{align*}
Here $Td_*$, $L_*$ and $c_*$ are respectively the Todd homology class, the $L$-homology class and the Chern homology class, i.e., the Poincar\'e dual of the corresponding characteristic class.

In particular, the constant coefficients and the top degree coefficients of them are respectively the following:
\begin{align*}
 & \chi^0(C_1\times \cdots   \times C_s \times S_1\times \cdots \times S_t) = \prod_{i=1}^s \chi^a(C_i) \cdot \prod_{j=1}^t \chi^a(S_j), \\
& \chi^{s+2t}(C_1\times \cdots   \times C_s \times S_1\times \cdots \times S_t) = \prod_{i=1}^s\frac{- \chi(C_i)}{2} \cdot \prod_{j=1}^t \frac{\sigma (S_j) + \chi(S_j) - 2\chi^a(S_j)}{2}
\end{align*}
\begin{align*}
 & T_{y*}^0(C_1\times \cdots   \times C_s \times S_1\times \cdots \times S_t) = \prod_{i=1}^s Td_*(C_i) \times  \prod_{j=1}^t Td_*(S_j), \\
& T_{y*}^{s+2t}(C_1\times \cdots  \times C_s \times S_1\times \cdots \times S_t) 
= \prod_{i=1}^s\frac{L_* (C_i) - c_*(C_i)}{2} \times  \prod_{j=1}^t \frac{L_*(S_j) + c_*(S_j) - 2Td_*(S_j)}{2}
\end{align*}
Here, as to the $T_{y*}^0$ and  $T_{y*}^{s+2t}$, the product $\prod_{i=1}^s Td_*(C_i)$ means the cross product and the same for the other products.
\end{thm}

\begin{rem} As far as the degree zero part and the top degree part are concerned, for the product of any varieties we can get the following formulae:
\begin{enumerate}
\item $ \displaystyle \chi^0(X_1 \times  \cdots \times X_n) = \prod_{i=1}^n \chi^a(X_i),$ \quad 
$ \displaystyle  \chi^{\sum_i \op{dim} X_i}(X_1 \times  \cdots  \times X_n) = \prod_{i=1}^n \chi^{\op{dim} X_i}(X_i),$
\item $ \displaystyle  T_{y*}^0(X_1 \times  \cdots \times X_n) = \prod_{i=1}^n Td_*(X_i),$ \quad 
$ \displaystyle  T_{y*}^{\sum_i \op{dim} X_i}(X_1 \times  \cdots \times X_n) = \prod_{i=1}^nT_{y*}^{\op{dim} X_i}(X_i).$
\end{enumerate}
In (2) the product $\prod_{i=1}^n $ is the cross product as in the above corollary.
\end{rem}

In the case when $n\geqq 3$, we need $n+1 \geqq 4$ points, thus the three points $-1, 0, 1$ are not enough. 



\begin{thm}\label{Vande}
Let $X$ be a compact nonsingular variety of dimension $n$. 
Let $a_0 =0, a_1=1, a_2 = -1, a_3, \cdots, a_n$ be mutually distinct numbers. Then we have 
$$\chi_y(X) = \chi^0(X)+ \chi^1(X)y + \chi^2(X)y^2+\cdots +\chi^n(X)y^n,$$
where
\begin{eqnarray*}
\left( 
\begin{array}{c}
\chi^0(X) \\
\chi^1(X) \\
\chi^2(X)  \\
\chi^3(X) \\
\vdots\\
\chi^{n-1}(X) \\
\chi^{n}(X)  \\
\end{array}
\right )
=
\left( 
\begin{array}{ccccccc}
1  & 0 & 0 & 0 \cdots \cdots & 0 \\
1  & 1 & 1 & 1 \cdots  \cdots & 1 \\
1  & -1 & (-1)^2 & (-1)^3 \cdots  \cdots & (-1)^n\\
1  & a_3 & a_3^2 & a_3^3 \cdots  \cdots & a_3^n\\
\vdots &&\ddots &&\vdots \\
1  & a_{n-1} & a_{n-1}^2 & a_{n-1}^3 \cdots \cdots & a_{n-1}^n\\
1  & a_n & a_n^2 & a_n^3 \cdots \cdots & a_n^n\\
\end{array} 
\right)^{\LARGE -1}
\left( 
\begin{array}{c}
\chi^a(X)\\
\sigma(X)\\
\chi(X) \\
\chi_{a_3}(X)\\
\vdots\\
\chi_{a_{n-1}}(X)\\
\chi_{a_n}(X)\\
\end{array}
\right)\\
\end{eqnarray*}
\end{thm}

\begin{ex} Let $n=3$ and $a_3 = 2$. (Here we note that the signature $\sigma(X) = 0$ since $\op{dim}_{\mathbb R}X = 6 \not \equiv 0 \op{mod} 4.$)
\begin{align*}
\chi_y(X) = \chi^a(X) + & \frac{-2\chi(X) -3\chi^a(X) - \chi_2(X)}{6} \, y \\
&+ \frac{\chi(X) - 2\chi^a(X)}{2} \, y^2 + \frac{\chi(X) + 3\chi^a(X) +\chi_2(X)}{6} \, y^3.
\end{align*}
\end{ex}

\begin{rem} Using the inversion formula in Remark \ref{inversion} (4), for example, if $n = 2k$, it suffices to consider the Vandermonde
$V_{n+1} \left (0,1,-1, 2, 3, \cdots, k, 2^{-1} ,3^{-1} \cdots, k^{-1} \right )$
and the special values $\chi_0(X)=\chi^a(X), \chi_1(X)=\sigma(X), \chi_{-1}(X)=\chi(X), \chi_2(X)$, $\chi_3(X)$, $\cdots$ and $\chi_k(X)$.
 \begin{eqnarray*}
\left( \begin{array}{c}
\chi^0(X) \\
\chi^1(X) \\
\chi^2(X) \\
\chi^3(X) \\
\chi^4(X) \\
\vdots\\
\chi^{k+1}(X) \\
\chi^{k+2}(X) \\
\chi^{k+3}(X) \\
\vdots\\
\chi^{2k}(X) \\
\end{array}
\right ) 
=
\left( 
\begin{array}{ccccccc}
1  & 0 & 0 & 0 & \cdots \cdots   & 0 \\
1  & 1 & 1 & 1 & \cdots  \cdots   & 1 \\
1  & -1 & 1 & -1 & \cdots  \cdots  & 1\\
1  & 2 & 2^2 & 2^3 & \cdots  \cdots  & 2^{2k}\\
1  & 3 & 3^2 & 3^3 & \cdots  \cdots  & 3^{2k}\\
\vdots &&\ddots &&\vdots \\
1  & k & k^2 & k^3 & \cdots  \cdots  & k^{2k}\\
1  & 2^{-1} & 2^{-2} & 2^{-3} & \cdots  \cdots  & 2^{-2k}\\
1  & 3^{-1} & 3^{-2} & 3^{-3} & \cdots  \cdots  & 3^{-2k}\\
\vdots &&\ddots &&\vdots \\
1  & k^{-1} & k^{-2} & k^{-3} & \cdots \cdots   & k^{-2k}\\
\end{array} 
\right)^{-1}
\left( 
\begin{array}{c}
\chi_0(X) \\
\chi_1(X) \\
\chi_{-1}(X) \\
\chi_2(X) \\
\chi_3(X) \\
\vdots\\
\chi_k(X) \\
2^{-2k}\chi_2(X) \\
3^{-2k}\chi_3(X) \\
\vdots\\
k^{-2k}\chi_k(X) \\
\end{array}
\right)\\
\end{eqnarray*}
\end{rem}

\begin{qu} Let $i$ be any integer greater than $1$. Can one express $\chi_i(X)$ and $\chi_{-i}(X)$ in terms of some other known invariants such as the Euler-Poincar\'e characteristic $\chi(X)$, the arithmetic genus $\chi^a(X)$, the signature $\sigma(X)$ and so on? 
\end{qu}

Theorem \ref{Vande}  also holds  for the motivic Hirzebruch class $T_{y*}(X) = T_{y*}([X \xrightarrow {\op{id}_X} X])$ for any possibly singular variety $X$ of dimension $n$:

\begin{thm} \label{Vande2} Let $X$ be a possibly singular variety of dimension $n$. Let $a_0 =0, a_1=1, a_2 = -1, a_3, \cdots, a_n$ be mutually distinct numbers. Then we have the following formula
$$T_{y*}(X) = T_*^0(X)+ T_*^1(X)y + T_*^2(X)y^2+\cdots +T_*^n(X)y^n,$$
where
\begin{eqnarray*}
\left( 
\begin{array}{c}
T_*^0(X) \\
T_*^1(X) \\
T_*^2(X)  \\
T_*^3(X) \\
\vdots\\
T_*^{n-1}(X)\\
T_*^n(X) \\
\end{array}
\right )
=
\left( 
\begin{array}{ccccccc}
1  & 0 & 0 & 0 \cdots \cdots & 0 \\
1  & 1 & 1 & 1 \cdots  \cdots & 1 \\
1  & -1 & (-1)^2 & (-1)^3 \cdots  \cdots & (-1)^n\\
1  & a_3 & a_3^2 & a_3^3 \cdots  \cdots & a_3^n\\
\vdots &&\ddots &&\vdots \\
1  & a_{n-1} & a_{n-1}^2 & a_{n-1}^3 \cdots \cdots & a_{n-1}^n\\
1  & a_n & a_n^2 & a_n^3 \cdots \cdots & a_n^n\\
\end{array} 
\right)^{\LARGE -1}
\left( 
\begin{array}{c}
T_{0*}(X)\\
T_{1*}(X)\\
c_*(X)\otimes \mathbb Q\\
T_{a_3*}(X)\\
\vdots\\
T_{a_{n-1}*}(X)\\
T_{a_{n}*}(X)\\\
\end{array}
\right)\\
\end{eqnarray*}
\end{thm}

We get the following corollary from \cite{BSY1} and \cite{MS}:
\begin{cor} Let the situation be as in Theorem \ref{Vande2}.
\begin{enumerate}
\item If $X$ is a toric variety, then $T_{0*}(X)$ can be replaced by Baum-Fulton-MacPherson's Todd class $td_*(X)$.
\item If $X$ is a simplicial projective toric variety, then $T_{0*}(X)$ can be replaced by Baum-Fulton-MacPherson's Todd class $td_*(X)$ and furthermore $T_{1*}(X)$ can be replaced by Capell-Shaneson's homology $L$-class $L_*(X)$.
\end{enumerate}
\end{cor}
\section{Derived Hirzebruch $\chi_y$-genus and derived motivic Hirzebruch class}

As stated above, the first motivation is trying to get a general formula of the Hirzebruch $\chi_y$-genus as well as the the motivic Hirzebruch class $T_{y*}(X)$, using the special values of them at $y=-1, 0, 1$ and some other points, the secondary motivation of the present paper is as follows: 

The Euler-Poincar\'e characteristic $\chi(X)$ is the alternating sum of the Betti numbers, i.e., $\chi(X) = \sum_i (-1)^i \op{dim} H_i(X; \mathbb R).$ If we use the Poincar\'e polynomial $P_X(t):= \sum_i \op{dim} H_i(X; \mathbb R)\, t^i$, then we have
$$\chi(X) = P_X(-1).$$
Namely, the Euler-Poincar\'e characteristic $\chi(X)$ is the constant term $a_0= P_X(-1)$ of the Taylor expansion of the Poincar\'e polynomial $P_X(t)$ at $t = -1$:
$$P_X(t) = a_0 + a_1 (t+1) + a_2(t+1)^2 + \cdots + a_k (t+1)^k + \cdots a_n (t+1)^n,$$
where $n$ is the degree of the Poincar\'e polynomial. More precisely $a_k = \frac {P^{(k)}_X(-1)}{k!}.$ The coefficient $a_1$ is called the secondary Euler characteristic, denoted by $\chi^{(1)}(X)$ and the other coefficients $a_k$ shall be called  the $k$-th higher Euler characteristic and denoted by $\chi^{(k)}(X)$ (cf. N. Ramachandran's recent paper \cite{Ramachandran} and \cite{GN}). We shall understand the Chern-Schwartz-MacPherson class $c_*(X)$ to be a homology class version of the the Euler-Poincar\'e characteristic, then a very 
natural question is the following:
\begin{qu} Is there a homology class version of the $k$-th higher Euler characteristic $\chi^{(k)}(X)$?
\end{qu}
To be able to answer this question, if we could have a certain reasonable ``Poincar\'e polynomial" for a compact variety $X$ 
$$\frak P_X(t) \in H_*(X; \mathbb R)[t] $$ 
such that
\begin{enumerate}
\item $\frak P_X(t) = c_*^{SM}(X) +  \frak a_1 (t+1) +  \frak a_2(t+1)^2 + \cdots +  \frak a_k (t+1)^k + \cdots  \frak a_n (t+1)^n,$
\item $\int_X \frak P_X(t)  = P_X(t)$, i.e., $\int_X c_*^{SM}(X) = \chi(X)$ and $\int_X \frak a_k = a_k,$
\end{enumerate}
then $\frak a_k$ would be such a homology class version of the $k$-th higher Euler characteristic $\chi^{(k)}(X)$.

At the moment we do not know if there would be such a ``Poincar\'e polynomial" $\frak P_X(t) 
$, although we have MacPherson's Chern class transformation $c_*:F(X) \to H_*(X)$. However, in the case of the Hirzebruch $\chi_y$-genus $\chi_y(X)$, which is a polynomial such that the special value of $\chi_y(X)$ at $y=-1$, i.e. $\chi_{-1}(X)$ is the Euler-Poincar\'e characteristic, we have the Taylor expansion of $\chi_y(X)$ at $y=-1$
$$\chi_y(X) = \chi_{-1}(X) + a_1(X)(y+1) + a_2(X) (y+1)^2 + a_3(X)(y+1)^3 +\cdots .$$
It turns out that many people (e.g., see \cite{LW}, \cite{NR}, \cite{Salamon}, \cite{Hir}, \cite{Li1}, \cite{Li2} etc.) have studied the first few coefficients $a_1(X), a_2(X), a_3(X), a_4(X)$. Thus, motivated by this fact, we also consider such coefficients for the motivic Hirzebruch classes, which we call \emph{derived motivic Hirzebruch classes}.\\

First we consider the following general situation. Given a polynomial of degree $n$:
$$f(y) = b_0 + b_1y + b_2y^2 + \cdots + b_py^p + \cdots + b_ny^n,$$
we consider the Taylor expansion of $f(y)$ at $y=\alp$:
$$f(y) = a_0 + a_1(y-\alp) + a_2(y-\alp)^2 + \cdots + a_p(y-\alp)^p + \cdots + a_n(y-\alp)^n,$$
where each coefficient $a_p$ can be expressed as follows
$$a_p= \frac{f^{(p)}(\alp)}{p!}.$$
We can see the following relations between $a_i$'s and $b_j$'s:
\begin{equation}\label{relation}
a_p = \sum_{k=p}^n  {k \choose p }  b_k \alp^{k-p} \quad \text{or} \quad b_p = \sum_{k=p}^n  {k \choose p }  a_k (-\alp)^{k-p}
\end{equation}
Expressing them in the matrix forms, we have
 \begin{eqnarray*}
\left( 
\begin{array}{ccccccc}
1  & \alp & \alp^2  & \alp^3 & \cdots & \cdots  \cdots  &  \alp^n \\
0  & 1 & {2\choose 1}\alp & {3\choose 1}\alp^2& \cdots  & \cdots  \cdots  & {n\choose 1}\alp^{n-1} \\
0  & 0 & 1 & {3\choose 2}\alp& \cdots  & \cdots \cdots  & {n\choose 2}\alp^{n-2}\\
0  & 0 & 0 & 1 & \cdots  & \cdots  \cdots & {n\choose 3}\alp^{n-3}\\
\vdots & &\ddots &&&&\vdots \\
0  & 0 & 0 & 0 & 1 & {p+1 \choose p}\beta \cdots   & {n\choose p}\alp^{n-p}\\
\vdots & &\ddots &&&&\vdots \\
0 & 0 & 0 & 0 & 0 &\cdots  \cdots  1 & {n \choose n-1}\alp\\
0 & 0 & 0 & 0 & 0 & \cdots \cdots  \cdots \cdots  & 1\\
\end{array} 
\right)
\left( \begin{array}{c}
b_0\\
b_1\\
b_2\\
b_3\\
\vdots\\
b_p\\
\vdots\\
b_{n-1}\\
b_n\\
\end{array}
\right )
= 
\left( \begin{array}{c}
a_0\\
a_1\\
a_2\\
a_3\\
\vdots\\
a_p\\
\vdots\\
a_{n-1}\\
a_n\\
\end{array}
\right)\\
\end{eqnarray*}
 \begin{eqnarray*}
\left( 
\begin{array}{ccccccc}
1  & -\alp  & \alp^2  & -\alp^3 & \cdots & \cdots  \cdots  &  (-\alp)^n \\
0  & 1 & -{2\choose 1}\alp & {3\choose 1}\alp^2& \cdots  & \cdots  \cdots  & {n\choose 1}(-\alp)^{n-1} \\
0  & 0 & 1 & -{3\choose 2}\alp & \cdots  & \cdots \cdots  & {n\choose 2}(-\alp)^{n-2}\\
0  & 0 & 0 & 1 & \cdots  & \cdots  \cdots & {n\choose 3}(-\alp)^{n-3}\\
\vdots & &\ddots &&&&\vdots \\
0  & 0 & 0 & 0 & 1 & -{p+1 \choose p}\alp\cdots   & {n\choose p}(-\alp)^{n-p}\\
\vdots & &\ddots &&&&\vdots \\
0 & 0 & 0 & 0 & 0 &\cdots  \cdots  1 & -{n \choose n-1}\alp\\
0 & 0 & 0 & 0 & 0 & \cdots \cdots  \cdots \cdots  & 1\\
\end{array} 
\right)
\left( \begin{array}{c}
a_0\\
a_1\\
a_2\\
a_3\\
\vdots\\
a_p\\
\vdots\\
a_{n-1}\\
a_n\\
\end{array}
\right )
= 
\left( \begin{array}{c}
b_0\\
b_1\\
b_2\\
b_3\\
\vdots\\
b_p\\
\vdots\\
b_{n-1}\\
b_n\\
\end{array}
\right)\\
\end{eqnarray*}

In particular $a_0 = f(\alp)$. In the case when $f(y) = \chi_y(X)$ is the Hirzebruch $\chi_y$-genus, i.e., 
$$\chi_y(X) = \chi^0(X) + \chi^1(X) y + \chi^2(X) y^2 + \cdots + \chi^p(X)y^p + \cdots + \chi^n(X)y^n,$$ 
(thus each $b_i = \chi^i(X)$), for a compact nonsingular variety or a compact rational homology manifold $X$  
$\chi_{-1}(X) = \chi(X)$ and $\chi_1(X) = \sigma(X)$ are respectively the first constant term $a_0, d_0$ of the following Taylor expansion of $\chi_y(X)$ at $y=-1, y=1$:
\begin{enumerate}
\item $\chi_y(X) = \chi_{-1}(X) + a_1(X)(y+1) + a_2(X) (y+1)^2 + \cdots + a_n(X)(y+1)^n,$
\item $\chi_y(X) = \chi_1(X) + d_1(X)(y-1) + d_2(X) (y-1)^2 + \cdots + d_n(X)(y-1)^n.$
\end{enumerate}
As to the Taylor expansion of $\chi_y(X)$ at $y=-1$, many people have already studied first few terms
$a_1(X), a_2(X), a_3(X), a_4(X)$ (e.g., see \cite{LW}, \cite{NR}, \cite{Salamon}, \cite{Hir}, \cite{Li1, Li2} etc.): for a compact complex manifold $V$ of dimension $n$
\begin{enumerate}
\item $\displaystyle a_1(V) = -\frac{1}{2}n\chi(V)$
\item $\displaystyle a_2(V) = \frac{1}{12} \left [ \frac{1}{2}n(3n-5)\chi(V) + c_{n-1}c_1\right ] $
\item $\displaystyle a_3(V) = -\frac{1}{24}\left [ \frac{n(n-2)(n-3)}{2}\chi(V) + (n-2)c_1c_{n-1} \right]$
\item 
$\displaystyle a_4(V) = -\frac{1}{5760} \Bigl [ n(15n^3  - 150n^2 + 485n - 502) \chi(V) 
+4(15n^2 -85n  + 108)c_1c_{n-1} + 8(c_1^2 +3c_2)c_{n-2} - 8(c_1^3 -3c_1c_2 +3c_3)c_{n-3} \Bigr] $
\end{enumerate}
Here we denote the Chern numbers $\left (c_{i_1}(V) \cdots c_{i_j}(V) \right )\cap [V]$ simply by $c_{i_1} \cdots c_{i_j}$, where $i_1 + \cdots i_j = n = \op{dim} V$.

\begin{rem} In \cite{LW} A. S. Libgober and J. W. Wood compute all the coefficients $a_i(V)\, \, (1\leqq i \leqq n)$ in the case when $1\leqq n=\op{dim} V \leqq 6$.
We recently found O. Debarre's paper \cite{D}, in which he computes all the coefficients $a_i(V)$ for $1\leqq n \leqq 9$.
\end{rem}

It follows from the above formula (\ref{relation}) that obtaining a general formula for the coefficient $\chi^p(X)$ is equivalent to obtaining a general formula for the coefficient $a_p(X)$, and thus it seems to be quite hard (or almost impossible) to get a general formula for $\chi_0(X) = \chi^0(X) = \chi^a(X)$. Which would suggest that it would be quite hard to get a general explicit formula for $\chi_y(X)$ like some cases done in \S 4. The merit of considering the Taylor expansion of $\chi_y(X)$ at \underline{$y = -1$} is that one can compute or express concretely at least first few or several coefficients $a_i(X)$ unlike the case of the coefficient $\chi^i(X)$. This kind of thing is called ``\emph{$-1$-phenomena}" in Ping Li's recent works \cite{Li1, Li2}.\\


Motivated by the higher Euler-Poincar\'e characteristic, we will introduce the following:

\begin{defn}
$$ \chi_y^{(p)}(X) := \frac{1}{p!} \frac{d^p}{dy^p} \Bigl(\chi_y(X) \Bigr).$$
is called \emph{the $p$-th derived Hirzebruch $\chi_y$-genus of $X$}.

$$T_{y*}^{(p)}(X) := \frac{1}{p!} \frac{d^p}{dy^p} \Bigl(T_{y*}(X)  \Bigr)$$
is called \emph{the $p$-th derived motivic Hirzebruch class of $X$}.
\end{defn}

As to the motivic Hirzebruch class \emph{transformation} $T_{y*}: K_0(\mathcal V/X) \to H_*(X) \otimes \mathbb Q[y]$, we define the \emph{$p$-th derived motivic Hirzebruch class transformation} as the following composition:
$$T_{y*}^{(p)}:= \frac{d^p}{dy^p} \circ T_{y*}: K_0(\mathcal V/X) \xrightarrow {T_{y*}} H_*(X) \otimes \mathbb Q[y] \xrightarrow {\frac{d^p}{dy^p}} H_*(X) \otimes \mathbb Q[y].$$
The naturality of this transformation is clear because of the naturality of the differential 
$\mathcal D := \frac{d}{dy}$, i.e, the commutativity of the following diagram for a proper morphisn $f:X \to Y$:
$$\CD
H_*(X) \otimes \mathbb Q[y]  @> {\mathcal D} >> H_*(X) \otimes \mathbb Q[y] \\
@V f_* VV @VV f_* V\\
H_*(Y) \otimes \mathbb Q[y] @>>  {\mathcal D}  > H_*(Y) \otimes \mathbb Q[y] . \endCD
$$
Which implies that for any $p>0$ and $\mathcal D^p = \overbrace{\mathcal D \circ \mathcal D \circ \cdots \circ \mathcal D}^p = \frac{d^p}{dy^p}$ the following commutes:
$$\CD
H_*(X) \otimes \mathbb Q[y]  @> {\mathcal D^p} >> H_*(X) \otimes \mathbb Q[y] \\
@V f_* VV @VV f_* V\\
H_*(Y) \otimes \mathbb Q[y] @>>  {\mathcal D^p}  > H_*(Y) \otimes \mathbb Q[y] . \endCD
$$

\begin{rem} We note the following equalities:
\begin{enumerate}
\item $\chi_y^{(p)}(X) = T_{y*}^{(p)}([X \to pt])$, where $X$ is compact.
\item $T_{y*}^{(p)}(X) = T_{y*}^{(p)}([X \xrightarrow {\op{id}_X} X]).$
\end{enumerate}
\end{rem}

\begin{rem}
Note that if $X$ is a toric variety, $\chi_y^{(p)}(X)$ and $T_y^{(p)}(X)$ are explicitly calculated
in \cite[Theorem 1.1 and Formula (1.7)]{MS}; e.g. $(-1)^p\cdot \chi_y^{(p)}(X)$ is just the number of the $p$-dimensional torus orbits.
\end{rem}
The following formula follows from the commutativity of the motivic Hirzebruch class and the cross products, and taking derivatives of cross products is similar to taking derivatives of product of two functions: 
\begin{thm} For two varieties $X, Y$, we have
$$T_{y*}^{(p)}([V \to X] \times [W \to Y])) = \sum_{i=0}^p {p \choose i} T_{y*}^{(i)}([V \to X]) \times T_{y*}^{(p-i)}([W \to Y])).$$
\end{thm}

\begin{cor}\label{formula 1} For two varieties the following holds:
\begin{enumerate}
\item $\displaystyle \chi_y^{(p)}(X \times Y) = \sum_{i=0}^p {p \choose i} \chi_y^{(i)}(X) \cdot \chi_y^{(p-i)}(Y)$, where $X$ and $Y$ are compact.
\item $\displaystyle  T_{y*}^{(p)}(X \times Y) = \sum_{i=0}^p {p \choose i} T_{y*}^{(i)}(X) \times T_{y*}^{(p-i)}(Y).$
\end{enumerate}
\end{cor}

\begin{rem} Note that the formula (1) of Corollary \ref{formula 1} follows directly from the multiplicativity of $\chi_y$.
\end{rem}
$\chi_y: K_0(\mathcal V/X) \to \mathbb Q[y]$ is a group (in fact, ring) homomorphism, thus the Euler-Poincar\'e characteristic $\chi_{-1} = \chi: K_0(\mathcal V/X) \to \mathbb Z (\subset \mathbb Q)$ is a homomorphism for complex algebraic varieties. We recall that the Euler-Poincar\'e characteristic is the alternating sum of the Betti numbers $b_i(X)$. In fact, in \cite{MP} C. McCrory and A. Parusi\'nski proved that for real algebraic varieties the Betti number can be ``captured" as a group homomorphism and thus the Poincar\'e polynomial can be also ``captured" as a group homomorphism.
\begin{thm}[Clint McCrory and Adam Parusi\'nski] There is a unique group homomorphism
$$\beta_i: K_0(\mathcal V_{\mathbb R}) \to \mathbb Z$$
such that for a compact variety $X$, $\beta_i(X) = b_i(X)$ is the usual Betti number.
$$\mathcal P(-)(t) := \sum_i \beta_it^i: K_0(\mathcal V_{\mathbb R}) \to \mathbb Z[t]$$ 
is a unique group homomorphism such that for a compact variety $X$
$$\mathcal P(X)(t) = \sum_i \beta_i(X) t^i = \sum_i b_i(X) t^i = P_X(t)$$
 is the usual Poincar\'e polynomial.
\end{thm}

Using this theorem, we can see that 
$$\chi_{(p)}(-):= \frac{d^p}{dt^p} \left (\mathcal P(-)(t)\right )|_{-1}:K_0(\mathcal V_{\mathbb R}) \to \mathbb Z $$
is a homomorphism version of Ramachandran's $p$-th higher Euler characteristic.

At the moment a natural transformation $?:K_0(\mathcal V_{\mathbb R}/X) \to \mathbb H_*(X)\otimes \mathbb Z[t]$ has not been constructed or found yet, thus we do not have a natural transformation version of 
Ramachandran's $p$-th higher Euler characteristic either, and the above $p$-th derived motivic Hirzebruch class transformation is the only one which is available.\\\\

\noindent
{\bf Acknowledgement.} We would like to thank the referee for many 
valuable comments and suggestions.\\

\end{document}